\documentclass[12pt,leqno,draft]{article}

\newtheorem{theorem}{Theorem}
\newtheorem{lemma}[theorem]{Lemma}
\newtheorem{proposition}[theorem]{Proposition}
\newtheorem{definition}[theorem]{Definition}
\newtheorem{corollary}[theorem]{Corollary}

\newcommand{\begintheorem}{\addtocounter{equation}{1}\begin{theorem}}
\newcommand{\beginlemma}{\addtocounter{equation}{1}\begin{lemma}}
\newcommand{\beginproposition}{\addtocounter{equation}{1}\begin{proposition}}
\newcommand{\begindefinition}{\addtocounter{equation}{1}\begin{definition}}
\newcommand{\begincorollary}{\addtocounter{equation}{1}\begin{corollary}}

\begin{document}

\title{Notes on normed algebras}

\author{Stephen William Semmes	\\
	Rice University		\\
	Houston, Texas}

\date{}

\maketitle

	All vector spaces and so forth here will be defined over the
complex numbers.  If $z = x + i \, y$ is a complex number, where $x$,
$y$ are real numbers, then the complex conjugate of $z$ is denoted
$\overline{z}$ and defined to be $x - i \, y$.  The complex conjugate
of a sum or product of complex numbers is equal to the corresponding
sum or product of complex conjugates.  The modulus of a complex number
$z$ is the nonnegative real number $|z|$ such that $|z|^2$ is equal to
the product of $z$ and its complex conjugate.  Thus the modulus of a
product of complex numbers is equal to the product of their moduli,
and one can show that the modulus of a sum of two complex numbers is
less than or equal to the sum of the moduli of the complex numbers.

	By a \emph{finite-dimensional algebra} we mean a finite
dimensional complex vector space $\mathcal{A}$ equipped with a binary
operation which satisfies the usual associativity and distributivity
properties, and which has a nonzero multiplicative identity element
$e$.  In other words, $e \, x = x \, e = x$ for all $x \in
\mathcal{A}$.  Thus $\mathcal{A}$ should have positive dimension in
particular.  Notice that the multiplicative identity element $e$ is
unique.

	As a basic class of examples, let $V$ be a finite-dimensional
complex vector space with positive dimension, and consider
$\mathcal{L}(V)$, the space of linear mappings from $V$ to itself.
This is a vector space whose dimension is equal to the square of the
dimension of $V$.  It also becomes an algebra with respect to the
usual composition of linear transformations, with the identity
transformation $I$ on $V$, which sends every vector in $V$ to itself,
as the multiplicative identity element.  If $\mathcal{A}$ is any
finite-dimensional algebra, then we can identify $\mathcal{A}$ with a
subalgebra of $\mathcal{L}(\mathcal{A})$, the algebra of linear
transformations on $\mathcal{A}$ considered simply as a vector space.
Namely, each element $a$ of $\mathcal{A}$ can be identified with the
linear transformation $x \mapsto a \, x$ on $\mathcal{A}$.

	As another class of examples, let $X$ be any finite nonempty
set.  Consider the vector space of complex-valued functions on $X$.
This becomes an algebra with respect to pointwise multiplication of
functions.  The multiplicative identity element for this algebra is
the function which is equal to $1$ at each point.  Of course this
algebra is commutative.

	As a third class of examples, let $A$ be a finite semigroup
with identity element $\theta$.  Thus $A$ is a finite set, $\theta$ is
an element of $A$, and there is a binary operation on $A$ which
associates to each pair of elements $x$, $y$ of $A$ another element $x
\, y$.  This operation should be associative, so that $x (y \, z)$ and
$(x \, y) z$ should be the same for all $x$, $y$, $z$ in $A$, and it
should satisfy $\theta \, x = x \, \theta = x$ for all $x \in A$.  As
usual, $\theta$ is uniquely determined by this feature.

	Consider the vector space of complex-valued functions on $A$.
If $f_1$, $f_2$ are two such functions, then we can define their
convolution to be the function on $A$ given by
\begin{equation}
	(f_1 * f_2)(z) = \sum_{x \, y = z} f_1(x) \, f_2(y).
\end{equation}
More precisely, this sum is taken over all $x, y \in A$ such that $x
\, y = z$.  In this way the functions on $A$ becomes an algebra, using
convolution as the multiplication operation.  The multiplicative
identity element is provided by the function which is equal to $1$
at the identity element $\theta$ in $A$ and equal to $0$ at all other
elements of $A$.

	If $A$ happens to be a commutative semigroup, then the
corresponding convolution algebra will also be commutative.  For each
element $x \in A$ we can define $\delta_x$ to be the function on $A$
which is equal to $1$ at $x$ and to $0$ at other elements of $A$, and
in this way we can embedd $A$ into its own convolution algebra in such
a way the multiplication in $A$ corresponds exactly to convolution of
the corresponding functions on $A$.  These functions $\delta_x$, $x
\in A$, form a basis for the vector space of functions on $A$, and
hence the dimension of the convolution algebra of functions on $A$ is
equal to the number of elements of $A$.

	By a norm on a finite-dimensional vector space $V$ we mean a
nonnegative real-valued function $N$ on $V$ such that $N(v) = 0$ if
and only if $v = 0$, $N(v + w) \le N(v) + N(w)$ for all $v, w \in V$,
and $N(\alpha \, v) = |\alpha| \, N(v)$ for all complex numbers
$\alpha$ and $v \in V$.  If $\mathcal{A}$ is a finite-dimensional
algebra and $\|\cdot \|$ is a norm on $\mathcal{A}$ as a vector space,
then we say that $(\mathcal{A}, \|\cdot \|)$ is a normed algebra if
also $\|x \, y \| \le \|x\| \, \|y\|$ for all $x, y \in \mathcal{A}$
and the multiplicative identity element $e \in \mathcal{A}$ has norm
equal to $1$.  For instance, if $V$ is a finite-dimensional complex
vector space with finite dimension and $N$ is a norm on $V$, then the
algebra $\mathcal{L}(V)$ of linear transformations on $V$ becomes a
normed algebra with repsect to the operator norm, which associates to
a linear transformation $T$ on $V$ the maximum of $N(T(v))$ over all
$v \in V$ with $N(v) = 1$.  If $X$ is any finite nonempty set, then
the algebra of functions on $X$ becomes a normed algebra when one uses
the norm which assigns to a complex-valued function $f$ on $X$ the
maximum of $|f(x)|$ over $x \in X$.  If $A$ is a finite semigroup,
then the convolution algebra becomes a normed algebra with respect to
the norm which assigns to a function $f$ on $A$ the sum of $|f(x)|$
over $x \in A$.

	Let $V$ be a finite-dimensional complex vector space of
positive dimension, and let $\mathcal{L}(V)$ denote the algebra of
linear transformations on $V$.  One can say that a linear mapping $T$
on $V$ is invertible if it is a one-to-one mapping of $V$ onto itself,
in which case the inverse of $T$ as a mapping on $V$ is also linear.
By well-known results in linear algebra $T$ is invertible if it is a
one-to-one mapping of $V$ into itself or if it maps $V$ onto itself,
because $V$ is finite-dimensional.  In algebraic terms $T$ is
invertible if there is a linear mapping $R$ on $V$ such that $R \, T$
and $T \, R$ are equal to the identity mapping on $V$.  Again because
$V$ has finite dimension, if either $R \, T$ or $T \, R$ is equal to
the identity mapping, then so is the other.

	If $T$ is a linear mapping on $V$, $I$ is the identity mapping
on $V$, and $\lambda$ is a complex number, then we get a new linear
mapping $\lambda \, I - T$.  The determinant of $\lambda \, I - T$ is
a complex number which is a polynomial in $\lambda$ of degree equal to
the dimension $n$ of $V$, and indeed the coefficient of $\lambda^n$ in
this polynomial is equal to $1$.  If $p(z)$ is any polynomial in $z$,
then we can define $p(T)$ in the usual manner, namely as the same
linear combination of powers of $T$ and the identity as we have of
powers of $z$ and $1$ in $p(z)$.  The celebrated theorem of Cayley and
Hamilton states that if $p(\lambda)$ is the polynomial given by taking
the determinant of $\lambda \, I - T$, then $p(T) = 0$.

	Of course $T$ is invertible if and only if the determinant of
$T$ is different from $0$.  If $n = 1$, then $T$ can be identified
with a complex number, and this is the same as saying that that
complex number is not equal to $0$.  When $n \ge 1$, the fact that
$p(T) = 0$ when $p(\lambda)$ is the characteristic polynomial equal to
the determinant of $\lambda I - T$ implies that if the determinant of
$T$ is different from $0$, so that the constant term of $p(\lambda)$
is different from $0$, then the inverse of $T$ can be expressed as a
polynomial of $T$ of degree $n - 1$.

	Now let $\mathcal{A}$ be a finite-dimensional algebra with
multiplicative identity element $e$.  If $x$ is an element of
$\mathcal{A}$, then we say that $x$ is invertible in $\mathcal{A}$ if
there is an element $y$ of $\mathcal{A}$ such that $y \, x = x \, y =
e$.  We can be a bit more precise and say that $x$ is left invertible
if there is an element $y_1$ of $\mathcal{A}$ such that $y_1 \, x =
e$, and that $x$ is right invertible if there is an element $y_2$ of
$\mathcal{A}$ such that $x \, y_2 = e$.  If $x$ is both left and right
invertible, with left and right inverses $R_1$, $R_2$, respectively,
then it is easy to see that $y_1 = y_2$ and $x$ is invertible.  In
particular, if $x$ is invertible, then the inverse of $x$ is unique,
and the inverse of $x$ is denoted $x^{-1}$.

	Suppose that $\mathcal{A}$ is in fact a subalgebra of
$\mathcal{L}(V)$ for some finite-dimensional vector space $V$ of
positive dimension.  As before, this can always be arranged up to
isomorphic equivalence.  If $T$ is an element of $\mathcal{A}$ which
is invertible as an element of $\mathcal{A}$, then of course $T$ is
invertible as an element of $\mathcal{L}(V)$.  Conversely, if $T$ is
invertible as an element of $\mathcal{L}(V)$, then the inverse of $T$
can be expressed as a polynomial in $T$, which therefore is an element
of $\mathcal{A}$.  In particular, it follows that $T$ is invertible if
$T$ is either left or right invertible.

	Fix a finite-dimensional algebra $\mathcal{A}$ with
multiplicative identity element $e$, and let $x$ be an element of
$\mathcal{T}$.  The spectrum of $x$ is defined to be the set of
complex numbers $\lambda$ such that $\lambda \, e - x$ does not have
an inverse in $\mathcal{A}$.  By embedding $\mathcal{A}$ into
$\mathcal{L}(V)$ for some finite-dimensional vector space $V$ we get
that the spectrum of $x$ can be described as the set of zeros of a
nonconstant polynomial on the complex numbers, and thus that the
spectrum of $x$ is a finite nonempty set of complex numbers.

	Assume that $(\mathcal{A}, \|\cdot \|)$ is a
finite-dimensional normed algebra.  Let $x$ be an element of
$\mathcal{A}$ and let $\lambda$ be a complex number such that
\begin{equation}
	|\lambda| > \|x\|.
\end{equation}
If $\mathcal{A}$ is a subalgebra of $\mathcal{L}(V)$ for some
finite-dimensional vector space $V$ and $\|\cdot \|$ is the operator
norm of a linear operator on $V$ with respect to some norm $N(\cdot )$
on $V$, then it is easy to see that $\lambda \, I - x$ has trivial
kernel as a linear operator on $V$ and hence is invertible.  In
general one can show that $\lambda \, e - x = \lambda (e -
\lambda^{-1} \, x)$ is invertible by summing the series $\sum_{j =
0}^\infty \lambda^{-j} \, x^j$.

\end{document}